\documentstyle[amssymb,10pt]{amsart}

\newtheorem{proposition}{Proposition}[section]

\newcounter{remark}[section]
\renewcommand{\theremark}{\arabic{section}.\arabic{remark}}
\newenvironment{remark}%
{\refstepcounter{remark}\trivlist \item[\hskip
    \labelsep {\bf Remark \theremark}]}%
{\endtrivlist}

\begin{document}

\newcommand{\nc}{\newcommand}
\nc{\pa}{\partial}
\nc{\cA}{{\cal A}}\nc{\cB}{{\cal B}}\nc{\cC}{{\cal C}}
\nc{\cE}{{\cal E}}\nc{\cG}{{\cal G}}\nc{\cH}{{\cal H}}\nc{\cZ}{{\cal
Z}}\nc{\cW}{{\cal W}} 
\nc{\cX}{{\cal X}}
\nc{\cY}{{\cal Y}}
\nc{\cR}{{\cal R}}\nc{\Ad}{\on{Ad}}
\nc{\id}{{\rm id}}
\nc{\ad}{\on{ad}}\nc{\Der}{\on{Der}}\nc{\End}{\on{End}}
\nc{\Imm}{\on{Im}}
\nc{\de}{\delta}\nc{\si}{\sigma}\nc{\on}{\operatorname}
\nc{\al}{\alpha}
\nc{\C}{{\Bbb C}}
\nc{\B}{{\cal B}}
\nc{\z}{{\Bbb Z}}
\nc{\la}{{\lambda}}
\nc{\wt}{\widetilde}
\nc{\G}{\wt{\g}}
\nc{\g}{{\frak g}}
\nc{\h}{{\frak h}}
\nc{\n}{{\frak n}}
\nc{\ep}{\epsilon}
\nc{\barf}{\overline{f}}
\nc{\barr}{\overline{R}}
\nc{\no}{\nonumber}
\newcommand{\uq}{{U_q(\widehat{sl}_2)}}
\newcommand{\baruq}{{\overline{U}_q(\widehat{sl}_2)}}
\newcommand{\ve}{{\varepsilon}}
\newcommand{\otimesdot}{{\stackrel{\textstyle.}{\otimes}}}
\newcommand{\x}{\xi}
\nc{\wz}{\textstyle{w\over z}}
\nc{\zw}{\textstyle{z\over w}}
\nc{\zab}{\textstyle{z_1\over z_2}}
\nc{\zba}{\textstyle{z_2\over z_1}}
\title[Singular $R$-matrices]
{Singular $R$-matrices and Drinfeld's comultiplication}
\author{R. Kedem}
\address{Department of Mathematics, University of California,
Berkeley, CA 94720}
%\date{October 1996}

\begin{abstract}
We compute the $R$-matrix which intertwines two dimensional evaluation
representations with Drinfeld comultiplication for $\uq$
\cite{Dr2}. This $R$-matrix contains terms proportional to the
$\delta$-function.  We construct the algebra $A(R)$ \cite{RSTS}
generated by the elements of the matrices $L^\pm(z)$ with relations
determined by $R$.  In the category of highest weight representations
there is a Hopf algebra isomorphism between $A(R)$ and an extension
$\baruq$ of Drinfeld's algebra.
\end{abstract}
\maketitle
\section{Introduction}

This note contains two results concerning the algebra $\uq$ with
Drinfeld comultiplication \cite{Dr2}. This algebra is presented in
terms of current generators, and as an algebra, it is isomorphic
\cite{Be} to the Drinfeld-Jimbo $\uq$ \cite{Dr,J}. However it has a
different comultiplication and therefore a different Hopf algebra
structure.

The first result presented here is a direct calculation of the
$R$-matrix which acts as an intertwiner of two dimensional evaluation
representations. The second is the construction of the $L$-matrix
algebra $A(R)$, presented in terms of relations between generators
$L^\pm(z)$, following the construction in \cite{RSTS}, and the
isomorphism between $A(R)$ and Drinfeld's algebra.

The quantum $R$-matrix which intertwines evaluation representations
with Drinfeld comultiplication  is completely determined, up to a
scalar factor, from the evaluation representation of Drinfeld's
generators and the action of the comultiplication on these generators,
from the equation
$$
\widetilde{R}(\zw) (\pi_z\otimes\pi_{w})\Delta(x) 
= (\pi_z\otimes\pi_{w})\Delta'(x) \widetilde{R}(\zw),\qquad x\in\uq
$$
(see section 2.3 for notation). These relations are considered here as
relations between formal power series. It is shown that this implies
that the $R$-martrix contains terms proportional to the
$\delta$-function.

This matrix is then used to construct the algebra $A(R)$, generated by
elements of the triangular matrices $L^\pm(z)$, with relations
determined by the $R$-matrix, of the form
\begin{eqnarray*}
 R({\zw}) L_1^\pm(z) L_2^\pm(w) 
&=&L_2^\pm(w) L_1^\pm(z)R(\zw),\\
R({\zw q^{-c}} )  L_1^+(z) L_2^- (w)& =&
L_2^- (w)L_1^+(z)R({\zw q^{c}}) 
\end{eqnarray*}
(see section 3 for the precise definition), in a
construction following the method of \cite{FRT,RSTS}.  This algebra has a
natural Hopf algebra structure, and the intertwining relation of the
$R$-matrix follows easily from the Yang-Baxter relation for $R$. We
construct the Hopf algebra isomorphism between $A(R)$ and an
extended version of Drinfeld's algebra, $\baruq$.

In the final section the Poisson limit $q\rightarrow1$ of the quantum
algebra $A(R)$ is computed. The relations in this Poisson algebra
are determined by the classical $r$-matrix.

It is straightforward to generalize the construction presented here to
the case of $U_q(\widehat{sl}_n)$ \cite{Ke}. However most of the new
features of the construction are already seen in the case $n=2$. 

\section{Drinfeld's ``new realization'' of $\uq$}

\subsection{Generators and relations}

Consider the Drinfeld realization of the algebra $\uq$ \cite{Dr2}
generated by
\begin{equation}\label{gen}
q^{\pm c},\quad \phi^\pm_{\pm n},\quad \xi^\pm_m, \qquad n\in
\z_{\geq0},\quad m\in\z ,
\end{equation}
with $q^{\pm c}$ central and relations expressed in terms of the
generating series
\begin{equation}\label{gs}
\xi^\pm(z) = \sum_{n\in\z} \xi_n^\pm z^{n} ,\qquad
\phi^\pm(z) = \sum_{n\geq0} \phi^\pm_{\pm n} z^{\pm n}
\end{equation}
as follows.
Let ${\cal B}=\C[[h]]$ be ring formal power
series in $h$, where $q=e^{h/2}$. Let
\begin{equation}
 g(z) =\frac{q^2 z - 1}{z-q^2} \quad \in{\cal B}[[z]]
\end{equation}
be an element of ${\cal B}$ with coefficients formal power series in
$z$. (Throughout this paper the notation $f(z)$ will be used to denote
functions in ${\cal B}[[z]]$, $f(z^{-1})\in {\cal B}[[z^{-1}]]$,
etc..)

The relations between the generating series (\ref{gs}) are the
following formal power series relations:
\begin{equation}\label{uq}
\begin{split}
&\phi^+(0)\phi^-(0) = \phi^-(0)\phi^+(0)=1,\\
&\left[\phi^\pm(z) ,\phi^\pm(w)\right]=0,\\
&\phi^+(z) \phi^-(w)= 
\frac{g(\zw q^{-c})}{g(\zw q^{c})}\phi^-(w)\phi^+(z)\\
&\phi^\pm(z)\xi^\pm(w) = g\left(\left(\zw\right)^{\pm1}\right)
 \xi^\pm(w)\phi^\pm(z) \\
&\phi^\pm(z)\xi^\mp(w) = 
g\left(\left(\zw q^{c}\right)^{\pm1}\right)^{-1} \xi^\mp(w)\phi^\pm(z) \\
&(z-q^{\pm2}w)\xi^\pm(z)\xi^\pm(w) =(q^{\pm2}z-w) \xi^\pm(w)\xi^\pm(z)\\
&\left[\xi^+(z),\xi^-(w)\right]= (q-q^{-1})\left( \delta(\zw q^{-c}) \phi^-(w)
-\delta(\zw q^{c}) \phi^+(z)\right).\\
\end{split}
\end{equation}
Here, the function $\delta(z)$ is a formal series,
$$\delta(z)=\sum_{n\in {\z}} z^n.$$ It can also be regarded as a
distribution which acts on the space of functions $f(z)$ regular at
$z=1$.

\subsection{Hopf algebra structure}

As an algebra, Drinfeld's algebra (\ref{uq}) is isomorphic to the
Drinfeld-Jimbo $\uq$ \cite{Dr,J,Be}. However it can be endowed with the
Hopf algebra structure of \cite{Dr2} which is different from that
found in \cite{Dr,J}. In this paper, the notation $\uq$ refers to
Drinfeld's algebra as a Hopf algebra.
\noindent{\bf Co-product:}
\begin{equation}
\label{cop}
\begin{split}
\Delta(q^c) &= q^c \otimes q^c,\\
\Delta(\xi^+(z)) &= \xi^+(z)\otimes 1 + \phi^+(z) 
\otimes \xi^+(z q^{c_1}), \\
\Delta(\xi^-(z)) &= 1 \otimes \xi^-(z) + \xi^-(zq^{c_2}) \otimes \phi^-(z),
 \\
\Delta(\phi^+(z))&=\phi^+(z)\otimes\phi^+(zq^{c_1}),\\
\Delta(\phi^-(z))&=\phi^-(zq^{c_2})\otimes\phi^-(z).\\
\end{split}
\end{equation}
Here, $c_i$, $i=1,2$ is the value of the central charge acting on the $i$-th
factor in the  tensor product.

\noindent{\bf Co-unit:}
\begin{equation}
\ve(q^c)=1,~~~~ \ve(\phi^{\pm}(z))=1, ~~~~~ 
\ve(\xi^\pm(z))=0.
\end{equation}
\noindent{\bf Antipode:}
\begin{equation}
\begin{split}
&S(q^c)~=~q^{-c},\\
&S(\phi^\pm(z))~=~\phi^{\pm}(z q^{-c})^{-1},\\
&S(\xi^+(z))~=~-\phi^{+}(z q^{-c})^{-1}\xi^+(zq^{-c}),\\
&S(\xi^-(z))~=~-\xi^-(zq^{-c})\phi^{-}(z q^{-c})^{-1}.
\end{split}
\end{equation}
These satisfy the condition $m\circ
(S\otimes {\rm id})\circ\Delta=\ve=m\circ({\rm id}\otimes
S)\circ\Delta,$ where $m (x\otimes y)=xy$.

\begin{remark}\label{no} {\em Normal ordering}:
Let ${\cal O}$ be the category of highest weight representations,
where the positive modes $\phi^+_n, \xi^\pm_n, \phi^-_{-n}, n>0$ act
nilpotently to the right. The product of two generating functions in
$\uq$ is well defined when it is normally ordered, i.e. all positive modes
on the right and the negative on the left. Then in the category ${\cal
O}$ normally ordered products of generating functions act 
as Laurent series.

The relations above describe a Hopf algebra structure in the category
${\cal O}$. The antipode should have the property $\langle x
v,w\rangle = \langle v, S(x) w\rangle$, $x\in \uq, v\in {M}^*, w\in
{M}$, where $M\in {\cal O}$, $M^*\in {\cal O}^*$ the right dual to
$M$, and  ${\cal O}^*$ is
the category of lowest weight modules. Therefore acting by $S$ on a
normally ordered product of generating functions gives an
anti-normally ordered product which acts in ${\cal O}^*$ as Laurent
series. 
\end{remark}

\subsection{The $R$-matrix}

Let $V=\End(\C^2)$ be the two dimensional representation of $sl_2$ with basis
$v_1,v_2$.  Consider the
two-dimensional evaluation representation of $\uq$, $V_z= {V} 
\otimes_\C
\C[z,z^{-1}]$.  The map 
$$
\pi_z : \uq |_{c=0} \rightarrow {V}[z,z^{-1}]$$ is defined by
\cite{DI}
\begin{eqnarray}\label{evala}
&\xi^\pm(w) \mapsto(q-q^{-1})~ \delta(\frac{w}{z})\ \sigma^\pm
\quad&\in
V[z,z^{-1}][[w, w^{-1}]]  \no\\
& \phi^+(w)\mapsto
\left(\begin{array}{cc}d(q^2\wz)&0\\0&d(\wz)^{-1}\end{array}\right)
&\in V [z^{-1}][[{w}]],\no\\
&\phi^-(w)\mapsto
\left(\begin{array}{cc}d(\zw)^{-1}&0\\0&d(q^2\zw)\end{array}\right)
&\in V [z][[{w}^{-1}]],
\end{eqnarray}
where
\begin{equation}
\label{dpm}
d(z)= \frac{1-z}{q-q^{-1}z}
\end{equation}
and $\sigma^\pm$ are the Pauli matrices. Here, the notation
$\C[z][[w]]$ indicates formal power series in $w$ with coefficients
polynomials in $z$, etc.

\begin{proposition}\label{Rmat}
There is a unique, up to a scalar multiple, $4\times4$ matrix
$\widetilde{R}(z)$ with diagonal entries in $\B[[z]]$ and off diagonal
elements in $\B[[z,z^{-1}]]$ which satisfies the intertwining relation
\begin{equation}
\label{intertwiner}
\widetilde{R}(z/z') (\pi_z\otimes\pi_{z'})\Delta(x) 
= (\pi_z\otimes\pi_{z'})\Delta'(x) \widetilde{R}(z/z'),~~ x\in\uq.
\end{equation}
%Here, we consider $(\pi_{z}\otimes\pi_{z'})\Delta(x)$ as a Laurent
%polynomial in $z/z'$ acting on formal power series by multiplication.
This unique solution is
$\widetilde{R}(z)=R_{12}^{-1}(z)$, where
\begin{equation}
\label{R}
R_{12}(z) = f(z)
\left(\begin{array}{cccc}
1 & & & \\
& d(z) & 0 & \\
&\gamma(z)\delta(z) & d(q^2 z) &  \\
& & & 1  \end{array} \right)~,
\end{equation}
where $\gamma(z)= (q-q^{-1})d(z)=\frac{(q-q^{-1})(1-z)}{q-q^{-1}z}$.
\end{proposition}
\begin{remark}
Here, $(1-z)^n\delta(z)$, $n\in \z_{\geq0}$ can be thought of as a
distribution acting on the space of functions with poles of order less
than or equal to $n$.
\end{remark}

The proof of proposition (\ref{Rmat}) is provided in the appendix.

The matrix $R$ will be used in the next section to define the algebra
$A(R)$. The scalar factor $f(z)$ is uniquely determined from the
requirement that the quantum determinant is central (see next
section), which gives
\begin{equation}\label{f}
f(z) = q^{1/2}(1-z) \widetilde{g}(z)^2 \quad \in \B[[z]],
\end{equation}
where
$$
\widetilde{g}(z) = \frac{(q^4 z;q^4)_\infty}{(q^2 z;q^4)_\infty},
\qquad (z,q)_n=\prod_{j=0}^{n}(1-z q^j).
$$
The function $f(z)$ in (\ref{f}) is the unique solution to the
difference equation
\begin{equation}\label{fdiff}
f(z)f(q^2 z)=d(q^2z)^{-1},\qquad f(z) \in \B[[z]].
\end{equation}

With this choice of $f(z)$, the matrix $R$ is equal to the evaluation of the
universal $R$-matrix \cite{R}.

\subsection{Extension of Drinfeld's algebra}
For what follows it will be useful to describe a slightly extended
version of Drinfeld's algebra, $\baruq$ \cite{FR}.  Consider the
generators $\alpha^\pm_{\pm n}$,
$$
\alpha^\pm(z) = \sum_{n\geq0} \alpha^\pm_{\pm n} z^{\pm n},
$$
such that
\begin{equation}\label{phial}
\phi^{\pm}(z) = \alpha^{\pm}(z)^{-1}\alpha^{\pm}(q^2 z)^{-1}.
\end{equation}
The algebra $\baruq$ is generated by $\xi^\pm(z), \alpha^\pm(z), q^c$
with relations as in (\ref{uq}) and
\begin{equation}
\begin{split}
&\left[\alpha^\pm(z),\alpha^\pm(w)\right]=0, \\
&\alpha^+(z)\xi^+(w)= d(\zw)^{-1}\xi^+(w)\alpha^+(z) ,\\
&\alpha^-(z)\xi^-(w)= d(q^2\wz)^{-1}\xi^-(w)\alpha^-(z) ,\\
&\alpha^+(z)\xi^-(w)= d(\zw q^c)\xi^-(w)\alpha^+(z) ,\\
&\alpha^-(z)\xi^+(w)= d(\wz q^{-c+2})\xi^+(w)\alpha^-(z) ,\\
&\alpha^+(z)\alpha^-(w)=\frac{f(zq^c/w)}{f(zq^{-c}/w)}\alpha^-(w)\alpha^+(z).
\end{split}
\end{equation}

The Hopf algebra structure extends to $\baruq$:
\begin{equation}
\begin{split}
\Delta(\alpha^+(z))&=\alpha^+(z)\otimes\alpha^+(zq^{c_1}), \\
\Delta(\alpha^-(z))&=\alpha^-(zq^{c_2})\otimes\alpha^-(z),\\
S(\alpha^\pm(z))&=\alpha^{\pm}(z q^{-c})^{-1},\\
\ve(\alpha^\pm(z))&=1.
\end{split}
\end{equation}

The two dimensional evaluation representation of $\alpha^\pm(z)$ is
obtained by solving the difference equation (\ref{phial}). The result is
\begin{equation}
\label{evalb}
\begin{split}
\pi_z\alpha^+(w)&=f(\wz)\left(\begin{array}{cc} 1 & 0 \\ 0 & d(\wz)
\end{array} \right), \\
\pi_z\alpha^-(w)&=\frac{1}{f(\zw)}
\left(\begin{array}{cc} 1 & 0 \\ 0 & d(q^2\zw)^{-1}
\end{array} \right) .
\end{split}
\end{equation}

\section{The quantum current algebra $A(R)$}

We define the Hopf algebra $A(R)$ using the same methods introduced in
\cite{RSTS}. 

Consider the algebra generated by the coefficients of $L^\pm_{ij}(z)$,
with $i,j=1,2$,  $L^\pm_{ii}(z)\in A(R)[[z^{\pm1}]]$, and the off
diagonal elements in $ A(R)[[z,z^{-1}]]$. The matrix $L^+$ ($L^-$) is
lower (upper) triangular.

The determining relations are
\begin{equation}
\label{AR}
\begin{split}
& R_{12}({\zw}) L_1^\pm(z) L_2^\pm(w) 
=L_2^\pm(w) L_1^\pm(z)R_{12}(\zw),\label{pp}\\
&R_{12}({\zw q^{-c}} )  L_1^+(z) L_2^- (w) =
L_2^- (w)L_1^+(z)R_{12}({\zw q^{c}}) .
\end{split}
\end{equation}

Here the matrix $R(z)$ is the same as in (\ref{R}) with $f(z)$ as in
(\ref{f}). We choose $f(z)$ by requiring that the
quantum determinants 
$${\cal D}^\pm(z)=L^\pm_{11}(z)L_{22}^\pm(z q^{-2})$$ are central.
The choice of the shift $q^{-2}$ in the definition of ${\cal D}^\pm(z)$
is motivated by the identification of the evaluation of $L^+(z)$ with
$R(z)$ (see below).

Finally, the algebra $A(R)$ is defined as the quotient of the algebra
generated by $L_{ij}^\pm(z)$ with relations (\ref{AR}) by the relation
${\cal D}^\pm(z)=1$.

\begin{remark}
This construction is motivated by the following
considerations. Suppose there exists a universal $R$-matrix for
$\baruq$ such that $L^+(z)=(\id\otimes\pi_z){\cal R}$,
$L^-(z)=(\id\otimes\pi_z){\cal R}_{21}^{-1}$ and $R(z/w) =
(\pi_z\otimes\pi_w ) {\cal R}$ \cite{R}. Then the relations (\ref{AR})
follow from the Yang-Baxter relation for ${\cal R}$. 
\end{remark}

The algebra $A(R)$ is a Hopf algebra with comultiplication
\begin{equation}\label{lcop}
\Delta L^+(z) = L^+(z) \otimesdot L^+(z q^{c_1}) ,\quad
\Delta L^-(z)= L^-(zq^{c_2}) \otimesdot L^-(z)
\end{equation}
which is easily seen to be an algebra homomorphism of 
(\ref{AR}). (Here the notation $\otimesdot$ indicates
matrix multiplication in $V$ and tensor product in $A(R)$.)
The co-unit acts as $\ve(L^{\pm}) = {1},$ and the action of the
antipode is $$ S(L^\pm(z)) = L^\pm(zq^{-c})^{-1}.$$

\begin{remark}\label{inverse}
Since ${\cal D}^\pm(z)=1$,
$L^\pm_{ii}(z)^{-1}\in A(R)$.  The matrices $S(L^\pm(z))
=L^\pm(z)^{-1}$ are well defined in the category ${\cal O}^*$, as they
should be according to remark (\ref{no}.
\end{remark}

Define the map $$\mu: A(R) \rightarrow \baruq$$ as
\begin{eqnarray*}
L^+(z) &\mapsto&
 \left({\begin{array}{cc} \alpha^+(z) & 0 \\ \xi^+(z)\alpha^+(z) &
 \alpha^+(q^2 z)^{-1}
\end{array}}\right),\\
L^-(z) &\mapsto& \left({\begin{array}{cc} \alpha^-(z) & -\alpha^-(z)\xi^-(z) 
\\ 0 & \alpha^-(q^2 z)^{-1}
\end{array}}\right).
\end{eqnarray*}

\begin{proposition}
The map $\mu$ is a Hopf algebra isomorphism.
\end{proposition}
The proof is by direct calculation.

By using the evaluation representation for the Drinfeld generators, it
is easily shown that
$$
\pi_w L^{+}(z) = R(z/w),\qquad \pi_w L^{-}(z) = R_{21}^{-1}(w/z).
$$
(Here, we use the notation
$(\pi_w L^+(z)_{ij})_{kl}= R_{ij,kl}(z/w)$ where the action of $R$ on
the basis vectors of $V$ is 
$R(z) v_i \otimes v_j = \sum_{k,l=1}^2 R_{ij,lk}(z)v_k\otimes v_l$.

The relations (\ref{AR}) in the evaluation representation are
therefore equivalent
to the Yang-Baxter equation for $R$:
\begin{equation}\label{ybe}
R_{12}(z) R_{13}(wz) R_{23}(w) = R_{23}(w)R_{13}(wz)R_{12}(z).
\end{equation}

The intertwining relation (\ref{intertwiner}) with $x=L^\pm(w)_{ij}$
is again a simple consequence of the Yang-Baxter relations
(\ref{ybe}), with the co-product as in (\ref{lcop}).

Note that this Hopf algebra isomorphism provides a trivial proof that
$\Delta$ is an algebra homomorphism of Drinfeld's algebra ({\it c.f.}
\cite{DI}).

\section{The Poisson algebra limit}

The Poisson limit $h\rightarrow 0$ of Drinfeld's algebra is obtained by
keeping $p=q^c$ and the generators $\xi^\pm$ and $\phi^\pm$ or
$\alpha^\pm$ fixed. In this limit the algebra $\baruq$ becomes a
Poisson algebra with
$$
\{a,b\} = \lim_{h\rightarrow0}  \frac{a\cdot b - b \cdot
a}{h}, \quad a,b\in \baruq
$$
Explicitly,
$$\phi^\pm (z) = \alpha^\pm (z) ^{-2}$$ 
and
\begin{equation}
\begin{split}
%&\{\phi^\pm(z),\phi^\pm(w) \} = 0 ,\quad 
&\{\alpha^\pm(z),\alpha^\pm(w) \} = 0 ,\\
%&\{\phi^+(z),\phi^-(w) \} = \left(\lambda(\zw p) - \lambda(\zw
%p^{-1})\right)\phi^+(z)\phi^-(w),\\
%&\{\phi^\pm(z),\xi^\pm(w) \} = -\lambda((\zw)^{\pm1}) \phi^\pm(z)\xi^\pm(w),\\
%&\{\phi^\pm(z),\xi^\mp(w) \} = \lambda((\zw p)^{\pm1}) \phi^\pm(z)\xi^\mp(w)
%,\\
&\{\xi^+(z),\xi^-(w) \} = \left( \delta(\zw p^{-1}) \phi^-(w) -
\delta(\zw p) \phi^+(z)\right),\\
&\{\xi^\pm(z),\xi^\pm(w)\}= \pm\frac{1}{2}\left(\lambda(\wz)-
\lambda(\zw)\right)\xi^\pm(z)\xi^\pm(w),\\
&\{\alpha^+(z),\alpha^-(w) \} = \frac{1}{4}\left(\lambda(\zw p) - \lambda(\zw
p^{-1})\right)\alpha^+(z)\alpha^-(w),\\
&\{\alpha^\pm(z),\alpha^\pm(w) \} =
 -\frac{1}{2}\lambda((\zw)^{\pm1}) \alpha^\pm(z)\xi^\pm(w),\\
&\{\alpha^\pm(z),\xi^\mp(w) \} = \frac{1}{2}\lambda((\zw p)^{\pm1}) 
\alpha^\pm(z)\xi^\mp(w).\\
\end{split}
\end{equation}
Again, these are formal power series relations with
$$
\lambda (z) = \frac{1}{4}\ \frac{1+z}{1-z}.
$$

The algebra $A(R)$ in this limit is a Poisson algebra generated by the
elements ${\cal L}_{ij}^\pm(z)$
as above,
and the quantum determinant is the usual determinant. The
relations (\ref{AR}) become, in this limit, 
\begin{equation}
\begin{split}
\left\{ {\cal L}_1^\pm(z),{\cal L}^\pm_2(w) \right\} &=
\left[{\cal L}_1^\pm(z){\cal L}^\pm_2(w),r(\zw)\right],\\
\left\{ {\cal L}_1^+(z),{\cal L}^-_2(w) \right\} &=
{\cal L}_1^+(z){\cal L}^-_2(w)r(\zw p^{-1})-r(\zw p)
{\cal L}_1^+(z){\cal L}^-_2(w).
\end{split}
\end{equation}

The classical $r$-matrix is obtained from $R(z)$ by taking the
limit $h\rightarrow 0$ of the elements of $R$. We find
\begin{equation}
R(z) = 1 + h\ r(z) + {\cal O}(h^2).
\end{equation}
Here we consider the off
diagonal element of $R$ to be in ${\cal B}[[z,z^{-1}]]$. This means
that the coefficient in front of the $\delta$ function in $r$
is the coefficient of $h$ in the expansion of $(q-q^{-1})f(z)/d(z)$,
which is 1. Thus, the classical $r$-matrix is 
\begin{equation}
r(z) = \left(\begin{array}{cccc}
\lambda(z) & 0 & 0 & 0 \\
0 & -\lambda(z) & 0 & 0 \\
0 & \delta(z) & -\lambda(z) & 0 \\
0 & 0 & 0 & \lambda(z) \end{array}\right)\ .
\end{equation}

The isomorphism $\mu$ gives, in this limit,
\begin{equation}
{\cal L}^+ = \left(\begin{array}{ll} \alpha^+(z) & 0 \\
\alpha^+(z)x^+(z) & \alpha^+(z)^{-1}\end{array}\right),
\quad
{\cal L}^- = \left(\begin{array}{ll} \alpha^-(z) & \alpha^-(z)x^-(z) \\
0 & \alpha^-(z)^{-1}\end{array}\right).
\end{equation}
Since all the operators commute there is no problem with normal
ordering in this limit.

\vskip.3in

\noindent{\bf Acknowledgements.} The author is indebted to
N. Reshetikhin for motivating this construction and for numerous
helpful discussions.  Thanks also to
I. Grojnowski for useful discussions.  This work is partly supported
by NSF grant DMS-9296120.

\vskip.3in

\appendix
\section{The intertwining relation}

The following is a proof of Proposition (\ref{Rmat}).

Let 
$$\widetilde{R}(z) = \left(\begin{array}{cccc}
a_1(z) & 0 & 0 & 0 \\
0 & b_1(z) & \gamma_1 \delta(z) &  0 \\
0 & \gamma_2 \delta(z) & b_2(z) &  0 \\
0 & 0 & 0 & a_2(z) \end{array}
\right),$$
with the diagonal elements in $\C[[z]]$ and the
off-diagonal elements in $\C[[z,z^{-1}]]$. Then, up to a scalar
multiple, the unique solution to the intertwining relation
(\ref{intertwiner}) is $R_{12}^{-1}(z)$.

The proof consists of using the evaluation representation (\ref{evala})
to explicitly compute
$$
\widetilde{R}(z/z') (\pi_z\otimes\pi_{z'}) \Delta(\xi^\pm(w)) =
(\pi_z\otimes\pi_{z'} )\Delta'(\xi^\pm(w)) \widetilde{R}(z/z').
$$
Writing out these relations, it is immediately apparant that the
off-diagonal terms of the matrix $\widetilde{R}$ must be proportional
to the $\delta$ function if they are nonzero.
Cancelling an overall factor of $\delta(w/z)$\footnote{In the sense
the coefficient in front of $\delta(w/z)$, since it is independent of
$w$, must be zero.}
and setting $z'=1$, this
results in eight relations between formal power series in $z$ and $z^{-1}$:
\begin{eqnarray}
\label{1}
&& a_1(z) d(q^2 z^{-1}) = b_1(z) + d(q^2 z)\gamma_2(z)
\delta(z),\\
\label{2}
&& a_1(z)=\gamma_1(z)\delta(z) +d(q^2 z)b_2(z),\\
\label{3}
&& b_1(z)+\gamma_1(z)\delta(z)/d(z^{-1}) = a_2(z)/ d(z),\\
\label{4}
&& \gamma_2(z)\delta(z) + b_2(z) /d(z^{-1}) =a_2(z),\\
\label{5}
&& b_1(z) +\gamma_1(z)\delta(z)/d(z^{-1})=a_1(z)/ d(z),\\
\label{6}
&& \gamma_2(z)\delta(z)+b_2(z)/ d(z^{-1})=a_1(z),\\
\label{7}
&& a_2(z)d(q^2 z^{-1})=b_1(z)+d(q^2 z)\gamma_2(z)
\delta(z),\\
\label{8}
&& a_2(z)=\gamma_1(z)\delta(z)+b_2(z) d(q^2 z).
\end{eqnarray}
The first four relations come from considering $\Delta(xi^+)$ in
(\ref{intertwiner}) and the last four from $\Delta(\xi^-)$.  

From (\ref{2}), $\gamma_1(z)=0$ since all other terms are in $\C[[z]]$
and there is nothing to cancel terms in $\C[[z^{-1}]]$ coming from
$\delta(z)$. Therefore, from (\ref{3}), $b_1(z) =
\frac{a_2(z)}{d(z)}$, and from (\ref{5}), $b_1(z) =
\frac{a_1(z)}{d(z)}$, so $a_2=a_1$. From (\ref{8}), $a_1(z) = b_2(z)
d(q^2 z)$.

As formal power series,
$$ \frac{1}{1-z^{-1}} = \sum_{n\leq0} z^n =
\delta(z) - \frac{z}{1-z},$$
and therefore
$$
d^{-1}(z^{-1}) = d(q^2 z) + (q-q^{-1}) \delta(z).
$$
From (\ref{4}), $$\gamma_2(z) \delta(z) = (q^{-1}-q)b_2(z)\delta(z).$$

These relations determine $R$ only up to a scalar multiple. Let
$$a_1(z) = \frac{a(z)}{1-z}.$$ Then 
$$b_1(z) = \frac{a(z)}{(1-z) d(z)},\quad b_2(z) = \frac{ a(z)}{(1-z)
d(q^2 z)},\quad \gamma_2(z) = \frac{(q^{-1}-q)a(z)}{(1-z) d(q^2z)}.$$
Thus $\widetilde{R}$ is identified with the matrix $R^{-1}(z)$ 
\begin{equation}\label{r1}
R_{12}^{-1}(z) = q^{-1/2} \widetilde{g}(z)^{-2}
\left[\begin{array}{cccc}
\frac{1}{1-z} & & & \\
& \frac{q-q^{-1}z}{(1-z)^2} & 0 & \\
& \gamma'(z)\delta(z) & \frac{1}{q^{-1}-q z} &  \\
& & & \frac{1}{1-z}  \end{array} \right]~,
\end{equation}
with
$$\gamma'(z) = \frac{q^{-1}-q}{(1-z) d(q^2 z)}=\frac{q^{-1}-q}{q^{-1}-q z},$$
with $a(z)= q^{-1/2} \widetilde{g}(z)^{-2}$.

\end{document}